\documentclass[11pt]{amsart}
\usepackage{amssymb,latexsym,amsmath,graphicx,graphics,epic,eepic}

\theoremstyle{plain}
\newtheorem{theorem}{Theorem}
\numberwithin{theorem}{section}
\newtheorem{lemma}[theorem]{Lemma}
\newtheorem{proposition}[theorem]{Proposition}

\theoremstyle{definition}

\begin{document}

\title{On the Slicing Genus of Legendrian Knots}

\author{Hao Wu}

\address{Department of mathematics and Statistics\\
         Lederle Graduate Research Tower\\
         710 North Pleasant Street\\
         University of Massachusetts\\
         Amherst, MA 01003-9305\\
         USA}

\email{wu@math.umass.edu}

\begin{abstract}
We apply Heegaard-Floer homology theory to establish generalized
slicing Bennequin inequalities closely related to a recent result
of T. Mrowka and Y. Rollin proved using Seiberg-Witten monopoles.
\end{abstract}

\maketitle

\section{Introduction}

Let $\xi$ be an oriented $2$-plane distribution on an oriented
$3$-manifold $M$. (Unless otherwise specified, all $3$-manifolds
in this paper are closed, connected and oriented.) $\xi$ is said to be a
contact structure on $M$ if there is a $1$-form $\alpha$ on $M$ so
that $\xi=\ker\alpha$, $d\alpha|_{\xi}>0$ and $\alpha\wedge
d\alpha>0$. Such a $1$-form is called a contact form for $\xi$.
And $(M,\xi)$ is called a contact $3$-manifold. A knot $K$ in a
contact $3$-manifold $(M,\xi)$ is called a Legendrian knot if it's
tangent to $\xi$. (Unless otherwise specified, all the knots in
this paper are oriented.) $(M,\xi)$ is said to be overtwisted if
there is an embedded disk $D$ in $M$ s.t. $\partial D$ is
Legendrian, but $D$ is transverse to $\xi$ along $\partial D$. If
$(M,\xi)$ is not overtwisted, then it's called tight. For example,
the standard contact structure $\xi_{st}$ on $S^3$ given by the
complex tangencies of the unit $3$-sphere in $\mathbb{C}^2$ is
tight. Overtwisted contact structures are kind of "soft", and are
completely classified up to isotopy by the homotopy type of the
underlying $2$-plane distribution. (See \cite{E1}.) Tight contact
structures display more rigidity, and possess more interesting
properties.

There are two "classical" invariants, the Thurston-Bennequin
number $tb(K)$ and the rotation number $r(K)$, for a Legendrian
knot $K$ in $(S^3,\xi_{st})$. These are generalized to
null-homologous Legendrian knots in any contact $3$-manifold (c.f.
\cite{E}). Let $K$ be a null-homologous Legendrian knot in a contact $3$-manifold
$(M,\xi)$, and $\Sigma\subset M$ a Seifert surface of $K$. Let
$K'$ be a knot obtained by pushing $K$ slightly in the direction
of a vector field that is transverse to $\xi$ along $K$. Then the
Thurston-Bennequin number $tb(K,\Sigma)$ is defined to be the
intersection number $\#(K'\cap\Sigma)$. Let $u$ be the positive
unit tangent vector field of $K$. Then the rotation number
$r(K,\Sigma)$ is defined to be the pairing $\langle c_1(\xi,u),
[\Sigma] \rangle$, where $[\Sigma]\in H_2(M,K)$ is the relative
homology class represented by $\Sigma$. If we reverse the
orientation of $K$, then $tb(K,\Sigma)$ is unchanged, and
$r(K,\Sigma)$ changes sign. Note that 
$tb(K,\Sigma)$ and $r(K,\Sigma)$ depend on $\Sigma$ only through the
relative homology class $[\Sigma]$. If $M$ is a homology spere,
then $H_2(M,K)=\mathbb{Z}$, and $tb$, $r$ are independent of
$\Sigma$. In this case, we suppress $\Sigma$ from the notation.

In \cite{Ben}, D. Bennequin proved the following Bennequin
inequality:

For any Legendrian knot $K$ in $(S^3,\xi_{st})$,
\begin{equation}\label{b}
tb(K)+|r(K)|\leq 2g(K)-1,
\end{equation}
where $g(K)$ is the genus of $K$.

In \cite{E}, Y. Eliashberg generalized \eqref{b} to any tight
contact $3$-manifold, and get:

For any null-homologous Legendrian knot $K$ in a tight contact
$3$-manifold, and any Seifert surface $\Sigma$ of $K$,
\begin{equation}\label{gb}
tb(K,\Sigma)+|r(K,\Sigma)| \leq -\chi(\Sigma).
\end{equation}

In \cite{Ru}, L. Rudolph strengthened \eqref{b} to the slicing
Bennequin inequality:

For any Legendrian knot $K$ in $(S^3,\xi_{st})$,
\begin{equation}\label{sb}
tb(K)+|r(K)|\leq 2g_s(K)-1,
\end{equation}
where $g_s(K)$ is the slicing genus of $K$.

Let $W$ be an oriented $4$-manifold with connected boundary
$\partial W =M$, and $\xi$ a contact structure on $M$. Assume that
$K$ is a Legendrian knot in $(M,\xi)$, and $\Sigma$ is an embedded
surface in $W$ bounded by $K$. In \cite{MR}, T. Mrowka and Y. Rollin extended the definitions of $tb$ and $r$ to this situation. The following is their construction. Let $v$ be a
vector field on $M$ transverse to $\xi$. Extend $v$ to a vector
field on $W$, and denote by $\{\varphi_t\}$ the flow of this
extended vector field. For a small $\varepsilon>0$, let
$K'=\varphi_{\varepsilon}(K)$ and
$\Sigma'=\varphi_{\varepsilon}(\Sigma)$. Then the intersection
number $\#(\Sigma\cap\Sigma')$ is well defined. The
Thurston-Bennequin number is defined to be
$tb(K,\Sigma)=\#(\Sigma\cap\Sigma')$. Note that $tb(K,\Sigma)$
depends on $\Sigma$ only through the relative homology class
$[\Sigma]\in H_2(W,K)$, and, when $\Sigma\subset M$, this
definition coincide with the previous definition of $tb$. Assume
$\mathfrak{s}\in Spin^{\mathbb{C}}(W)$, and there is an
isomorphism $h:\mathfrak{s}|_M\rightarrow\mathfrak{t}_{\xi}$,
where $\mathfrak{t}_{\xi}$ is the canonical
$Spin^{\mathbb{C}}$-structure on $M$ associated to $\xi$. Choose a
complex structure on $\xi$. Then $\det(\mathfrak{t}_{\xi})$ is
canonically isomorphic to $\xi$, and $h$ induces an isomorphism
$\det(h):\det(\mathfrak{s})|_M\rightarrow\xi$. Let $u$ be the
positive unit tangent vector field of $K$. The rotation
number is defined to be $r(K,\Sigma,\mathfrak{s},h)=\langle
c_1(\det(\mathfrak{s}), \det(h)^{-1}(u)),[\Sigma]\rangle$. Note
that $r(K,\Sigma,\mathfrak{s},h)$ depends on $\Sigma$ only through
the relative homology class $[\Sigma]\in H_2(W,K)$, depends on the
pair $(\mathfrak{s},h)$ only through the isomorphism type of it in
$Spin^{\mathbb{C}}(W,\xi)$, and, again, when $\Sigma\subset M$,
$r$ is independent of $(\mathfrak{s},h)$ and coincide with the
previous definition of the rotation number. As before, under the
reversal of the orientation of $K$, $tb$ is unchanged, and $r$
changes sign. In the special case that there is a symplectic form
$\omega$ on $W$ such that $(W,\omega)$ is a weak symplectic
filling of $(M,\xi)$, i.e., $\omega|_{\xi}>0$, this symplectic
form $\omega$ determines a canonical $Spin^{\mathbb{C}}$-structure
$\mathfrak{s}_{\omega}$ on $W$ and a canonical isomorphism
$h_{\omega}:\mathfrak{s}_{\omega}|_M\rightarrow\mathfrak{t}_{\xi}$.
Write
$r(K,\Sigma,\omega)=r(K,\Sigma,\mathfrak{s}_{\omega},h_{\omega})$.

In \cite{MR}, T. Mrowka and Y. Rollin prove the following
generalized slicing Bennequin inequality using Seiberg-Witten
monopole invariants.

\begin{theorem}[\cite{MR}, Theorem A]\label{MR-main}
Let $W$ be an oriented $4$-manifold with connected boundary
$\partial W =M$, and $\xi$ a contact structure on $M$. Let $K$ be
a Legendrian knot in $(M,\xi)$, and $\Sigma$ an embedded surface
in $W$ bounded by $K$. Assume there is an element
$(\mathfrak{s},h)\in Spin^{\mathbb{C}}(W,\xi)$, such that
$SW(\mathfrak{s},h)\neq0$. Then
\begin{equation}\label{gsb-sw}
tb(K,\Sigma)+|r(K,\Sigma,\mathfrak{s},h)| \leq -\chi(\Sigma).
\end{equation}
Specially, when $(W,\omega)$ is a weak symplectic filling of
$(M,\xi)$ (c.f. \cite{KM}, Theorem 1.1), we have
\begin{equation}\label{gsb-weak-sw}
tb(K,\Sigma)+|r(K,\Sigma,\omega)| \leq -\chi(\Sigma).
\end{equation}
\end{theorem}

There are two approaches in the study of $3$-dimensional gauge
theory: the Seiberg-Witten-Floer approach by counting solutions to
the Seiberg-Witten equation; and the Heegaard-Floer approach by
counting holomorphic curves. Though the techniques used in these
two approaches are quite different, it is conjectured that these
give equivalent theories as their $4$-dimensional counterparts do.
In this paper, we use Heegaard-Floer homology to prove the
following generalizations of the slicing Bennequin inequality,
which further demonstrates the similarity between the two
theories.

\begin{theorem}\label{gsb-main}
Let $W$ be an oriented $4$-manifold with connected boundary
$\partial W =M$, $\xi$ a contact structure on $M$, and $K$ a
Legendrian knot in $(M,\xi)$.

(a) If there is a $Spin^{\mathbb{C}}$-structure $\mathfrak{s}$ on
$W$ with $F^+_{W\setminus B,\mathfrak{s}|_{W\setminus
B}}(c^+(\xi))\neq0$, where $B$ is an embedded $4$-ball in the
interior of $W$, then there is an isomorphism
$h:\mathfrak{s}|_M\rightarrow\mathfrak{t}_{\xi}$ such that, for
any embedded surface $\Sigma$ in $W$ bounded by $K$,
\begin{equation}\label{gsb-hf}
tb(K,\Sigma)+|r(K,\Sigma,\mathfrak{s},h)| \leq -\chi(\Sigma).
\end{equation}

(b) If $(W,\omega)$ is a weak symplectic filling of $(M,\xi)$,
then, for any embedded surface $\Sigma$ in $W$ bounded by $K$,
\begin{equation}\label{gsb-weak-hf}
tb(K,\Sigma)+|r(K,\Sigma,\omega)| \leq -\chi(\Sigma).
\end{equation}
\end{theorem}

\section{Heegaard-Floer Homology}

In this section, we review aspects of the Heegaard-Floer theory
necessary for the proof of Theorem \ref{gsb-main}.

\subsection{Heegaard-Floer homology}

In \cite{OS2}, P. Ozsv\'ath and Z. Szab\'o defined the
Heegaard-Floer homology groups of $3$-manifolds.
Given a connected oriented closed $3$-manifold $M$ and a
$Spin^{\mathbb{C}}$-structure $\mathfrak{t}$ on $M$, there are
four Heegaard-Floer homology groups associated to $M$:
$HF^\infty(M,\mathfrak{t})$, $HF^-(M,\mathfrak{t})$,
$HF^+(M,\mathfrak{t})$ and $\widehat{HF}(M,\mathfrak{t})$. The first three
are $\mathbb{Z}[U]$-modules, and the last one is a $\mathbb{Z}$-module. In this paper, we will mostly use
$HF^+(M,\mathfrak{t})$. Moreover, given a
$\mathbb{Z}[H^1(M)]$-module $\mathfrak{M}$, there is the notion of
$\mathfrak{M}$-twisted Heegaard-Floer homology
$\underline{HF}^+(M,\mathfrak{t};\mathfrak{M})$, which is a
$\mathbb{Z}[U]\otimes\mathbb{Z}[H^1(M)]$-module (c.f. \cite{OS3}).

If $\mathfrak{M}_1$ and $\mathfrak{M}_2$ are two
$\mathbb{Z}[H^1(M)]$-modules, and $\theta:\mathfrak{M}_1
\rightarrow \mathfrak{M}_2$ is a homomorphism, then $\theta$
naturally induces a homomorphism
\[
\Theta: \underline{HF}^+(M,\mathfrak{t};\mathfrak{M}_1)
\rightarrow \underline{HF}^+(M,\mathfrak{t};\mathfrak{M}_2).
\]

If we consider $\mathbb{Z}$ as a $\mathbb{Z}[H^1(M)]$-module, then
$\underline{HF}^+(M,\mathfrak{t};\mathbb{Z})$ is the (untwisted)
Heegaard-Floer homology $HF^+(M,\mathfrak{t})$ defined with the
appropriate coherent orientation system, and the $2^{b_1(M)}$
choices of $\mathbb{Z}[H^1(M)]$-module structures on $\mathbb{Z}$
correspond to the $2^{b_1(M)}$ coherent orientation systems on the
moduli spaces (c.f. \cite{OS2,OS3}).

In \cite{OS1}, P. Ozsv\'ath and Z. Szab\'o introduced the
Heegaard-Floer homology twisted by a $2$-form. More precisely,
consider the polynomial ring
\[
\mathbb{Z}[\mathbb{R}]=\{\sum_{i=1}^k c_iT^{s_i} ~|~
k\in\mathbb{Z}_{\geq0}, ~c_i\in\mathbb{Z}, ~s_i\in\mathbb{R}\}.
\]
Let $[\omega]\in H^2(M;\mathbb{R})$. The action $e^{[\nu]} \cdot
T^s = T^{s+ \int_M \nu\wedge\omega}$, where $[\nu]\in H^1(M)$,
gives $\mathbb{Z}[\mathbb{R}]$ a $\mathbb{Z}[H^1(M)]$-module
structure. Denote the module by
$\mathbb{Z}[\mathbb{R}]_{[\omega]}$. Then the Heegaard-Floer
homology of $M$ twisted by $[\omega]$ is defined to be
\[
\underline{HF}^+(M,\mathfrak{t};[\omega]) =
\underline{HF}^+(M,\mathfrak{t};\mathbb{Z}[\mathbb{R}]_{[\omega]}).
\]

In \cite{OS3}, P. Ozsv\'ath and Z. Szab\'o deduced the following
adjunction inequality:

\begin{theorem}[\cite{OS3}, Theorem 7.1]\label{ai}
Let $\Sigma$ be a close oriented surface embedded in a
$3$-manifold $M$ with $g(\Sigma)\leq1$, $\mathfrak{t}$ a $Spin^{\mathbb{C}}$-structure
on $M$, and $\mathfrak{M}$ a $\mathbb{Z}[H^1(M)]$-module. If
$\underline{HF}^+(M,\mathfrak{t};\mathfrak{M})\neq0$, then
\begin{equation}
|\langle c_1(\mathfrak{t}),[\Sigma]\rangle| \leq -\chi(\Sigma).
\end{equation}
\end{theorem}

Note that, although the adjunction inequality is only prove for
untwisted Heegaard-Floer homology in \cite{OS3}, the proof there
readily adapts to the twisted case.

\subsection{Homomorphisms induced by cobodisms}

Let $W$ be a cobodism from a $3$-manifold $M_1$
to another $3$-manifold $M_2$, and
$\mathfrak{s}$ a $Spin^{\mathbb{C}}$-structure on $W$. Then $W$
and $\mathfrak{s}$ induce a homomorphism
\[
F^+_{W,\mathfrak{s}}: HF^+(M_1,\mathfrak{s}|_{M_1}) \rightarrow
HF^+(M_2,\mathfrak{s}|_{M_2}).
\]

Let $\mathfrak{M}$ be a $\mathbb{Z}[H^1(M_1)]$-module, and
$\delta: H^1(\partial W) \rightarrow H^2(W,\partial W)$ the
connecting map in the long exact sequence of the pair $(W,\partial
W)$. Define
\[
\mathfrak{M}(W) = \mathfrak{M} \otimes_{\mathbb{Z}[H^1(M_1)]}
\mathbb{Z}[\delta H^1(\partial W)],
\]
where the action of $\mathbb{Z}[H^1(M_i)]$ on $\mathbb{Z}[\delta
H^1(\partial W)]$ is induced by $e^{[\nu]} \mapsto
e^{\delta([\nu])}$. Then $W$ and $\mathfrak{s}$ also induce a
homomorphism
\[
\underline{F}^+_{W,\mathfrak{s}}:
\underline{HF}^+(M_1,\mathfrak{s}|_{M_1};\mathfrak{M}) \rightarrow
\underline{HF}^+(M_2,\mathfrak{s}|_{M_2};\mathfrak{M}(W)).
\]
The definition of this homomorphism depends on some auxiliary
choices. So it's only define up to right action by units of
$\mathbb{Z}[H^1(M_1)]$ and left action by units of
$\mathbb{Z}[H^1(M_2)]$. Alternatively, we consider it as an
equivalence class of homomorphisms from
$\underline{HF}^+(M_1,\mathfrak{s}|_{M_1};\mathfrak{M})$ to
$\underline{HF}^+(M_2,\mathfrak{s}|_{M_2};\mathfrak{M}(W))$, and
denote this equivalence class by
$[\underline{F}^+_{W,\mathfrak{s}}]$.

Specially, for an $[\omega]\in H^2(W;\mathbb{R})$, let
$\mathfrak{M}=\mathbb{Z}[\mathbb{R}]_{[\omega|_{M_1}]}$. There is
a natural homomorphism $\theta: \mathfrak{M}(W) \rightarrow
\mathbb{Z}[\mathbb{R}]_{[\omega|_{M_2}]}$ induced by $T^s \otimes
e^{\delta([\nu])} \mapsto T^{s+\int_{M_2} \nu\wedge\omega}$, where
$[\nu]\in H^1(M_2)$. This map induces a homomorphism $\Theta$
between the Heegaard-Floer homologies of $M_2$ twisted by these
two $\mathbb{Z}[H^1(M_2)]$-modules. Composing it with
$\underline{F}^+_{W,\mathfrak{s}}$, we get a homomorphism
\[
\underline{F}^+_{W,\mathfrak{s};[\omega]} = \Theta \circ
\underline{F}^+_{W,\mathfrak{s}}:
\underline{HF}^+(M_1,\mathfrak{s}|_{M_1};[\omega|_{M_1}])
\rightarrow
\underline{HF}^+(M_2,\mathfrak{s}|_{M_2};[\omega|_{M_2}]).
\]
Again, $\underline{F}^+_{W,\mathfrak{s};[\omega]}$ is only defined
up to multiplication by $\pm T^s$, and is consider as a equivalent
class $[\underline{F}^+_{W,\mathfrak{s};[\omega]}]$ (c.f.
\cite{OS1}).

The homomorphisms defined here satisfy the following composition
laws:

\begin{theorem}[\cite{OS4}, Theorems 3.4, 3.9]\label{comp}
Let $W_1$ be a cobodism from a $3$-manifold
$M_1$ to another $3$-manifold $M_2$, and $W_1$
a cobodism from $M_2$ to a third $3$-manifold
$M_3$. Then $W=W_1\cup_{M_2}W_2$ is a cobodism from $M_1$ to
$M_3$. We have:

(a) For any $Spin^{\mathbb{C}}$-structures $\mathfrak{s}_i \in
Spin^{\mathbb{C}}(W_i)$, $i=1,2$,
\begin{equation}\label{untwisted-comp}
F^+_{W_2,\mathfrak{s}_2} \circ F^+_{W_1,\mathfrak{s}_1} =
\sum_{\{\mathfrak{s}\in Spin^{\mathbb{C}}(W) ~|~
\mathfrak{s}|_{W_i}\cong\mathfrak{s}_i\}} \pm
F^+_{W,\mathfrak{s}}.
\end{equation}

(b)  Let $\mathfrak{s} \in Spin^{\mathbb{C}}(W)$, and
$\mathfrak{s}_i = \mathfrak{s}|_{W_i}$. For any
$\mathbb{Z}[H^1(M_1)]$-module $\mathfrak{M}$, there are
representatives $\underline{F}^+_{W_1,\mathfrak{s}_1} \in
[\underline{F}^+_{W_1,\mathfrak{s}_1}]$ and
$\underline{F}^+_{W_2,\mathfrak{s}_2} \in
[\underline{F}^+_{W_2,\mathfrak{s}_2}]$ such that
\begin{equation}\label{twisted-comp}
[\underline{F}^+_{W,\mathfrak{s}}] =[\Pi \circ
\underline{F}^+_{W_1,\mathfrak{s}_1} \circ
\underline{F}^+_{W_2,\mathfrak{s}_2}],
\end{equation}
where $\Pi$ is induced by the natural homomorphism from
$\mathfrak{M}(W_1)(W_2)$ to $\mathfrak{M}(W)$.
\end{theorem}

Combine Theorem \ref{comp} and the blow-up formula (\cite{OS4},
Theorem 3.7), we have the following theorem.

\begin{theorem}[\cite{OS4}, Theorems 3.4, 3.7]\label{blow-up}
Let $W$ be a cobodism from a $3$-manifold $M_1$
to another $3$-manifold $M_2$, and
$\mathfrak{s}$ a $Spin^{\mathbb{C}}$-structure on $W$. Blow up an
interior point of $W$. We get a new cobodism $\widehat{W}$ from
$M_1$ to $M_2$. Let $\widehat{\mathfrak{s}}$ be the lift of
$\mathfrak{s}$ to $\widehat{W}$ with $\langle
c_1(\widehat{\mathfrak{s}}), [E]\rangle = -1$, where $E$ is the
exceptional sphere. Then
$F^+_{W,\mathfrak{s}}=F^+_{\widehat{W},\widehat{\mathfrak{s}}}$.
\end{theorem}

\subsection{The contact invariant}

In \cite{OS}, P. Ozsv\'ath and Z. Szab\'o defined the
Ozsv\'ath-Szab\'o invariants for contact $3$-manifolds. For each
contact $3$-manifold $(M,\xi)$, it is an element $c(\xi)$ of the
quotient $\widehat{HF}(-M,\mathfrak{t}_{\xi})/\{\pm1\}$, where
$\mathfrak{t}_{\xi}$ is the $Spin^{\mathbb{C}}$-structure
associated to $\xi$. Let $\iota:\widehat{HF}(-M)\rightarrow
HF^+(-M)$ be the natural map (c.f. \cite{OS2}). We set
$c^+(\xi)=\iota(c(\xi))\in HF^+(-M,\mathfrak{t}_{\xi})/\{\pm1\}$.
This version of the Ozsv\'ath-Szab\'o contact invariants is easier
to use for our purpose. Given a $\mathbb{Z}[H^1(M)]$-module
$\mathfrak{M}$, one can similarly define the twisted
Ozsv\'ath-Szab\'o contact invariant $c^+(\xi;\mathfrak{M}) \in
HF^+(-M,\mathfrak{t}_{\xi};\mathfrak{M})/\mathbb{Z}[H^1(M)]^{\times}$,
where $\mathbb{Z}[H^1(M)]^{\times}$ is the set of units of
$\mathbb{Z}[H^1(M)]$. Specially, if $[\omega]\in
H^2(M;\mathbb{R})$, then we have the $[\omega]$-twisted invariant
$c^+(\xi;[\omega]) \in HF^+(-M,\mathfrak{t}_{\xi};[\omega])/\{\pm
T^s|~s\in\mathbb{R}\}$ (c.f. \cite{OS1}). These contact invariants
vanish when $\xi$ is overtwisted. Following properties of the
Ozsv\'ath-Szab\'o contact invariants are needed for the proof of
of Theorem \ref{gsb-main}.

\begin{proposition}[\cite{Gh1}, Proposition 3.3]\label{Gh1-3.3}
Suppose that $(M',\xi')$ is obtained from $(M,\xi)$ by Legendrian
surgery on a Legendrian link. Then we have
$F^+_{W,\mathfrak{s}_0}(c^+(\xi'))=c^+(\xi)$, where $W$ is the
cobordism induced by the surgery and $\mathfrak{s}_0$ is the
canonical $Spin^\mathbb{C}$-structure associated to the symplectic
structure on $W$. Moreover, $F^+_{W,\mathfrak{s}}(c^+(\xi'))=0$
for any $Spin^\mathbb{C}$-structure $\mathfrak{s}$ on $W$ with
$\mathfrak{s}\ncong\mathfrak{s}_0$.
\end{proposition}

\begin{theorem}[\cite{OS1}, Theorem 4.2]\label{weak-fill-nv}
Let $(M,\xi)$ be a contact $3$-manifold with a weak symplectic
filling $(W,\omega)$. Let $B$ be an embedded $4$-ball in the
interior of $W$. Consider $W\setminus B$ as a cobodism from $-M$
to $-\partial B$. Then $\underline{F}^+_{W\setminus B,
\mathfrak{s}_{\omega}|_{W\setminus B};[\omega|_{W\setminus
B}]}(c^+(\xi;[\omega|_M]))\neq0$, where $\mathfrak{s}_{\omega}$ is
the $Spin^\mathbb{C}$-structure on $W$ associated to $\omega$.
\end{theorem}

\section{Generalized Slicing Bennequin Inequalities}

In this section, we adapt T. Mrowka and Y. Rollin's idea into the
Heegaard-Floer setting, and prove Theorem \ref{gsb-main}.

\begin{lemma}[\cite{MR}]\label{adjust}
Let $W$ be an oriented $4$-manifold with connected boundary
$\partial W=M$, $\xi$ a contact structure on $M$, and
$\mathfrak{s}$ a $Spin^\mathbb{C}$-structure on $W$ with an
isomorphism $h:\mathfrak{s}|_M\rightarrow\mathfrak{t}_{\xi}$.
Assume $K$ is a Legendrian knot in $(M,\xi)$, and $\Sigma \subset
W$ is an embedded surface bounded by $K$. Then there are a
Legendrian knot $K'$ in $(M,\xi)$ and an embedded surface $\Sigma'
\subset W$ bounded by $K'$, such that $tb(K',\Sigma')\geq1$,
$\chi(\Sigma')\leq-1$, and $tb(K',\Sigma') +
|r(K',\Sigma',\mathfrak{s},h)| + \chi(\Sigma') = tb(K,\Sigma) +
|r(K,\Sigma,\mathfrak{s},h)| + \chi(\Sigma)$.
\end{lemma}
\begin{proof}
Let $p$ be a point on $K$. There is a neighborhood $U$ of $p$ so that
$(U,\xi|_U)\cong(\mathbb{R}^3,\xi_0)$, where $\xi_0$ is the
standard contact structure on $\mathbb{R}^3$ defined by $dz-ydx$.
By the following Legendrian Reidemeister move, we create a pair of
cusps on the front projection of $K\cap U$ (c.f. \cite{FT}).

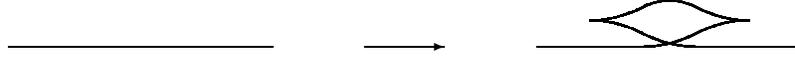
\begin{figure}[h]

\setlength{\unitlength}{1pt}

\begin{picture}(420,30)(-210,-10)

\linethickness{.5pt}


\put(-150,0){\line(1,0){100}}


\put(-15,0){\vector(1,0){30}}


\put(50,0){\line(1,0){40}}

\qbezier(90,0)(100,0)(110,5)

\qbezier(110,5)(120,10)(130,10)

\qbezier(130,10)(120,10)(110,15)

\qbezier(110,15)(100,20)(90,15)

\qbezier(70,10)(80,10)(90,15)

\qbezier(70,10)(80,10)(90,5)

\qbezier(110,0)(100,0)(90,5)

\put(110,0){\line(1,0){40}}

\end{picture}
\caption{Creating cusps}\label{reidemeister-move}
\end{figure}

Near a cusp, connect sum $K$ with a Legendrian righthand trefoil
knot $T_r$ in $U$ with $tb(T_r)=1$. We get a new Legendrian knot
$K_1$ and an embedded surface $\Sigma_1$ in $W$ bounded by $K_1$,
s.t., $tb(K_1,\Sigma_1)=tb(K,\Sigma)+1$,
$\chi(\Sigma_1)=\chi(\Sigma)-1$, and
$|r(K_1,\Sigma_1,\mathfrak{s},h)| = |r(K,\Sigma,\mathfrak{s},h)|$.
Repeat this process, we will find a $K'$ and a $\Sigma'$ with the
properties specified in the lemma.
\end{proof}

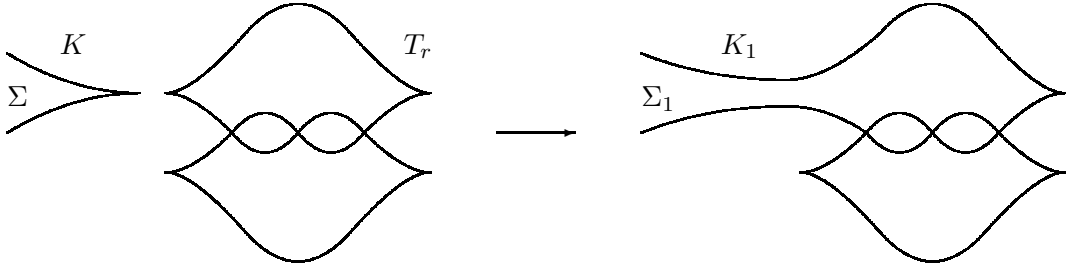
\begin{figure}[h]

\setlength{\unitlength}{1pt}

\begin{picture}(420,100)(-210,-50)

\linethickness{.5pt}


\qbezier(-200,30)(-175,15)(-150,15)

\qbezier(-200,0)(-175,15)(-150,15)

\put(-180,30){$K$}

\put(-200,10){$\Sigma$}

\qbezier(-140,15)(-130,15)(-110,37.5)

\qbezier(-110,37.5)(-90,60)(-70,37.5)

\qbezier(-70,37.5)(-50,15)(-40,15)

\qbezier(-140,15)(-130,15)(-115,0)

\qbezier(-115,0)(-102.5,-15)(-90,0)

\qbezier(-90,0)(-77.5,15)(-65,0)

\qbezier(-65,0)(-50,-15)(-40,-15)

\qbezier(-140,-15)(-130,-15)(-110,-37.5)

\qbezier(-110,-37.5)(-90,-60)(-70,-37.5)

\qbezier(-70,-37.5)(-50,-15)(-40,-15)

\qbezier(-140,-15)(-130,-15)(-115,0)

\qbezier(-115,0)(-102.5,15)(-90,0)

\qbezier(-90,0)(-77.5,-15)(-65,0)

\qbezier(-65,0)(-50,15)(-40,15)

\put(-50,30){$T_r$}


\put(-15,0){\vector(1,0){30}}


\qbezier(40,30)(65,20)(95,20)

\qbezier(40,0)(65,10)(95,10)

\qbezier(95,20)(110,20)(130,37.5)

\qbezier(130,37.5)(150,60)(170,37.5)

\qbezier(170,37.5)(190,15)(200,15)

\qbezier(95,10)(110,10)(125,0)

\qbezier(125,0)(137.5,-15)(150,0)

\qbezier(150,0)(162.5,15)(175,0)

\qbezier(175,0)(190,-15)(200,-15)

\qbezier(100,-15)(110,-15)(130,-37.5)

\qbezier(130,-37.5)(150,-60)(170,-37.5)

\qbezier(170,-37.5)(190,-15)(200,-15)

\qbezier(100,-15)(110,-15)(125,0)

\qbezier(125,0)(137.5,15)(150,0)

\qbezier(150,0)(162.5,-15)(175,0)

\qbezier(175,0)(190,15)(200,15)

\put(70,30){$K_1$}

\put(40,10){$\Sigma_1$}

\end{picture}
\caption{Connect summing with $T_r$}\label{connecting-sum}
\end{figure}

\vspace{.4cm}

\noindent\textit{Proof of Theorem \ref{gsb-main}.} By Lemma
\ref{adjust}, we only need prove the theorem for $K$ and $\Sigma$
with $tb(K,\Sigma)\geq1$ and $\chi(\Sigma)\leq-1$. We assume these
are true throughout the proof.

We prove part (a) first. Performing Legendrian surgery along $K$
gives a symplectic cobodism $(V,\omega')$ from $(M,\xi)$ to
another contact $3$-manifold $(M',\xi')$ (c.f. \cite{We,Wu}). By
Proposition \ref{Gh1-3.3},
$F^+_{V,\mathfrak{s}_{\omega'}}(c^+(\xi'))=c^+(\xi)$. Let
$\widetilde{W}=W\cup_{M}V$. Then, by Theorem \ref{comp}, we have
\[
\sum_{\{\widetilde{\mathfrak{s}}\in
Spin^{\mathbb{C}}(\widetilde{W}) ~|~
\widetilde{\mathfrak{s}}|_{W}\cong\mathfrak{s},
~\widetilde{\mathfrak{s}}|_{V}\cong\mathfrak{s}_{\omega'}\}} \pm
F^+_{W,\widetilde{\mathfrak{s}}}(c^+(\xi'))= F^+_{W,\mathfrak{s}}
\circ
F^+_{V,\mathfrak{s}_{\omega'}}(c^+(\xi'))=F^+_{W,\mathfrak{s}}(c^+(\xi))\neq
0.
\]
Thus, there is an $\widetilde{\mathfrak{s}}\in
Spin^{\mathbb{C}}(\widetilde{W}) $ with
$\widetilde{\mathfrak{s}}|_{W}\cong\mathfrak{s}$, and
$\widetilde{\mathfrak{s}}|_{V}\cong\mathfrak{s}_{\omega'}$, such
that
\[
F^+_{\widetilde{W}\setminus
B,\widetilde{\mathfrak{s}}|_{{W}\setminus B}}(c^+(\xi'))\neq0.
\]
Let $h_1:\widetilde{\mathfrak{s}}|_{W}\rightarrow\mathfrak{s}$,
$h_2:\widetilde{\mathfrak{s}}|_{V}\rightarrow\mathfrak{s}_{\omega'}$
be the above isomorphisms, and
$h_3:\mathfrak{s}_{\omega'}|_M\rightarrow\mathfrak{t}_{\xi}$ the
natural isomorphism. And define
$h:\mathfrak{s}|_M\rightarrow\mathfrak{t}_{\xi}$ by $h=h_3\circ
h_2 \circ h_1^{-1}$.

Capping off $\Sigma$ by the core of the $2$-handle, we get an
embedded closed surface $\widetilde{\Sigma}$ satisfying
$\chi(\widetilde{\Sigma})=\chi(\Sigma)+1\leq0$,
$[\widetilde{\Sigma}]\cdot[\widetilde{\Sigma}]=tb(K,\Sigma)-1\geq0$,
and $\langle
c_1(\widetilde{\mathfrak{s}}),[\widetilde{\Sigma}]\rangle=r(K,\Sigma,\mathfrak{s},h)$.
Next, blow up $tb(K,\Sigma)-1$ points on the core of the
$2$-handle, we get a new $4$-manifold $\widehat{W}$ with a natural
projection $\pi:\widehat{W}\rightarrow W$. Let
$\widehat{\mathfrak{s}}$ be the lift of $\widetilde{\mathfrak{s}}$
to $\widehat{W}$ whose evaluation on each exceptional sphere is
$-1$, and $\widehat{\Sigma}$ be the lift of $\widetilde{\Sigma}$
to $\widehat{W}$ obtained by removing the exceptional spheres from
$\pi^{-1}(\widetilde{\Sigma})$. Then
$\chi(\widehat{\Sigma})=\chi(\widetilde{\Sigma})=\chi(\Sigma)+1$,
$[\widehat{\Sigma}]\cdot[\widehat{\Sigma}]=0$, and $\langle
c_1(\widehat{\mathfrak{s}}),[\widehat{\Sigma}]\rangle =
r(K,\Sigma,\mathfrak{s},h)+tb(K,\Sigma)-1$. Also, by Theorem
\ref{blow-up},
\[
F^+_{\widehat{W}\setminus
\widehat{B},\widehat{\mathfrak{s}}|_{\widehat{W}\setminus
\widehat{B}}}(c^+(\xi')) = F^+_{\widetilde{W}\setminus
B,\widetilde{\mathfrak{s}}|_{\widetilde{W}\setminus
B}}(c^+(\xi'))\neq0,
\]
where $\widehat{B}\subset \widehat{W}$ is the pre-image of
$B\subset W\subset\widetilde{W}$ under $\pi$. Since
$[\widehat{\Sigma}]\cdot[\widehat{\Sigma}]=0$, there is a
neighborhood $U$ of $\widehat{\Sigma}$ in $\widehat{W}$
diffeomorphic to $\widehat{\Sigma}\times D^2$. Since the location
of $\widehat{B}$ does not affect the map
$F^+_{\widehat{W}\setminus
\widehat{B},\widehat{\mathfrak{s}}|_{\widehat{W}\setminus
\widehat{B}}}$, we assume that $\widehat{B}$ is in the interior of
$U$. Let $W_1=\widehat{W}\setminus U$, and $W_2=U\setminus
\widehat{B}$. Then, by Theorem \ref{comp}, there are maps
\[
\underline{F}^+_{W_1,\widehat{\mathfrak{s}}|_{W_1}}:
HF^+(-M',\mathfrak{t}_{\xi'}) \rightarrow
\underline{HF}^+(-\partial U, \widehat{\mathfrak{s}}|_{\partial
U};\mathbb{Z}(W_1)),
\]
\[
\underline{F}^+_{W_2,\widehat{\mathfrak{s}}|_{W_2}}:
\underline{HF}^+(-\partial U, \widehat{\mathfrak{s}}|_{\partial
U};\mathbb{Z}(W_1)) \rightarrow \underline{HF}^+(-\partial B,
\widehat{\mathfrak{s}}|_{\partial B};\mathbb{Z}(W_1)(W_2)),
\]
such that $F^+_{\widehat{W}\setminus
\pi^{-1}(B),\widehat{\mathfrak{s}}|_{\widehat{W}\setminus
\pi^{-1}(B)}} = \Theta \circ
\underline{F}^+_{W_2,\widehat{\mathfrak{s}}|_{W_2}} \circ
\underline{F}^+_{W_1,\widehat{\mathfrak{s}}|_{W_1}}$, where
\[
\Theta:\underline{HF}^+(-\partial B, \mathfrak{s}|_{\partial
B};\mathbb{Z}(W_1)(W_2))\rightarrow HF^+(-\partial B,
\mathfrak{s}|_{\partial B})
\]
is induced by the natural projection
$\theta:\mathbb{Z}(W_1)(W_2)\rightarrow\mathbb{Z}$. Specially,
this implies $\underline{HF}^+(-\partial U,
\widehat{\mathfrak{s}}|_{\partial U};\mathbb{Z}(W_1))\neq0$. Note
that $\partial U\cong \widehat{\Sigma}\times S^1$. Hence, by
Theorem \ref{ai}, we have
\[
\langle c_1(\widehat{\mathfrak{s}}),[\widehat{\Sigma}]\rangle \leq
-\chi(\widehat{\Sigma}),
\]
that is
\[
tb(K,\Sigma)+r(K,\Sigma,\mathfrak{s},h) \leq -\chi(\Sigma).
\]
Reverse the orientations of $K$ and $\Sigma$, and repeat the whole
arguement. We get
\[
tb(K,\Sigma)-r(K,\Sigma,\mathfrak{s},h) \leq -\chi(\Sigma).
\]
Thus,
\[
tb(K,\Sigma)+|r(K,\Sigma,\mathfrak{s},h)| \leq -\chi(\Sigma).
\]

Now we use twisted Heegaard-Floer homology to prove part (b).
Again, perform Legendrian surgery along $K$. This gives a new
contact $3$-manifold $(M',\xi')$ with a weak symplectic filling
$(\widetilde{W},\widetilde{\omega})$ (c.f. \cite{We,Wu}). Define
$\widehat{\Sigma}$ and $\widehat{W}$ as above, i.e., by capping
off $\Sigma$ with the core of the $2$-handle, and then blowing up
$tb(K,\Sigma)-1$ points on the core the of two handle. Let
$\widehat{\omega}$ be the blown-up symplectic form on
$\widehat{W}$. Then $(\widehat{W},\widehat{\omega})$ is also a
weak symplectic filling of $(M',\xi')$. Denote by
$\widehat{\mathfrak{s}}$ the canonical
$Spin^{\mathbb{C}}$-structure associated to $\widehat{\omega}$. We
have $\chi(\widehat{\Sigma})=\chi(\Sigma)+1$,
$[\widehat{\Sigma}]\cdot[\widehat{\Sigma}]=0$, and $\langle
c_1(\widehat{\mathfrak{s}}),[\widehat{\Sigma}]\rangle =
r(K,\Sigma,\omega)+tb(K,\Sigma)-1$. Also, by Theorem
\ref{weak-fill-nv}, we have $\underline{F}^+_{\widehat{W}\setminus
\widehat{B}, \widehat{\mathfrak{s}}|_{\widehat{W}\setminus
\widehat{B}};[\widehat{\omega}|_{\widehat{W}\setminus
\widehat{B}}]}(c^+(\xi';[\widehat{\omega}|_{M'}]))\neq0$. Again,
assume that $\widehat{B}$ is in the interior of $U$. Then, similar
to part (a), $\underline{F}^+_{\widehat{W}\setminus \widehat{B},
\widehat{\mathfrak{s}}|_{\widehat{W}\setminus
\widehat{B}};[\widehat{\omega}|_{\widehat{W}\setminus
\widehat{B}}]}$ factors through $\underline{HF}^+(-\partial U,
\widehat{\mathfrak{s}}|_{\partial
U};\mathbb{Z}[\mathbb{R}]_{[\widehat{\omega}|_{M'}]}(W_1))$, where
$W_1=\widehat{W}\setminus U$. So
\[
\underline{HF}^+(-\partial U, \widehat{\mathfrak{s}}|_{\partial
U};\mathbb{Z}[\mathbb{R}]_{[\widehat{\omega}|_{M'}]}(W_1))\neq0,
\]
and we apply Theorem \ref{ai} as above to prove part (b). \hfill
$\Box$

\end{document}